\documentclass[12pt]{article}

\usepackage[cmtip,arrow]{xy}
\usepackage{amsmath,amssymb,pb-diagram,pb-xy}%,pdfsync}
\def \a{{\mathfrak a}}
\def \al{\alpha}

\def \coker{{\rm coker}}
\def \coim{{\rm coim}}
\def \C{{\mathbb C}}

\def \CF{{\cal F}}
\def \CG{{\cal G}}
\def \CH{{\cal H}}
\def \CO{{\cal O}}

\def \F{{\mathbb F}}

\def \GL{{\rm GL}}
\def \Hom{{\rm Hom}}

\def \im{{\rm im}}

\def \Mod{{\rm Mod}}

\def \N{{\mathbb N}}
\def \P{{\mathbb P}}
\def \p{{\mathfrak p}}
\def \ph{\varphi}

\def \prf{{\bf Proof: }}
\def \q{{\mathfrak q}}
\def \Q{{\mathbb Q}}
\def \Quot{{\rm Quot}}
\def \qed{\ifmmode\eqno \square 
		\else\noproof\vskip 12pt plus 3pt minus 9pt \fi}
\def \noproof{{\unskip\nobreak\hfill\penalty50\hskip2em		\hbox{}%
     \nobreak\hfill $\square$\parfillskip=0pt%
     \finalhyphendemerits=0\par}}
\def \R{{\mathbb R}}

\def \setminus{\begin{picture}(18,10)\put(4,6)
                {\line(2,-1){10}}\end{picture}}

\def \spec{{\rm spec\,}}

\def \Z{{\mathbb Z}}
\def \({\left(}
\def \){\right)}
\def \={{\ =\ }}

\newcommand{\dotcup}{\ensuremath{\mathaccent\cdot\cup}}

\newtheorem{theorem}{Theorem}[section]

\newtheorem{lemma}[theorem]{Lemma}

\newtheorem{proposition}[theorem]{Proposition}

\begin{document}

\pagestyle{myheadings} \markright{REMARKS ON $F_1$-SCHEMES}

\title{$\F_1$-schemes and toric varieties}
\author{Anton Deitmar}
\date{}
\maketitle

%{\it At a time when so many scholars in the world are calculating, is it not
%desirable that some, who can, dream?}\ \footnote{Thom, Ren\'e: Structural Stability and Morphogenesis Addison Wesley 1989, p. 325}

$ $

{\bf MCC classification: 14A15}, 11G25, 14L32, 14M32


{\bf Author's address:} Math. Inst., Auf der Morgenstelle 10, 72076 T\"ubingen, Germany

{\bf Abstract:} In this paper it is shown that integral $\F_1$-schemes of finite type are essentially the same as toric varieties.
A description of the $\F_1$-zeta function in terms of toric geometry is given.
Etale morphisms and universal coverings are introduced.

%\tableofcontents

\section*{Introduction}
There are by now several attempts to make the theory of the field of one element $\F_1$ rigorous.
In \cite{ToenVaquie} the authors formalize the transition from rings to schemes on a categorial level and apply this machinery to the category of sets to obtain the category of $\F_1$-schemes as in \cite{F1}.
In \cite{Durov} and \cite{Haran} the authors extend the definition of rings in order to capture a structure that deserves to be called $\F_1$.
In \cite{F1} the author tried instead to fix the minimum properties any of these theories must share.
The current paper extends this line of thought.
%We show that schemes over $\F_1$ in the sense of \cite{F1}, which are of finite type, give rise to toric varieties after base change to $\C$.
%This is not true for the concepts of Durov \cite{Durov} or Haran \cite{Haran}.
We use terminology of \cite{F1} and \cite{zetaF1}.

In this paper, a ring will always be commutative with unit and a monoid will always be commutative.
An \emph{ideal} $\a$ of a monoid $A$ is a subset with $A\a\subset \a$.
A \emph{prime ideal} is an ideal $\p$ such that $S_\p=A\setminus\p$ is a submonoid of $A$.
For a prime ideal $\p$ let $A_\p=S_\p^{-1}A$ be the \emph{localization} at $\p$.
The \emph{spectrum} of a monoid $A$ is the set of all prime ideals with the obvious Zariski-topology (see \cite{F1}).
Similar to the theory of rings, one defines a structure sheaf $\CO_X$ on $X=\spec(A)$, and one defines a \emph{scheme over $\F_1$} to be a topological space together with a sheaf of monoids, locally isomorphic to spectra of monoids.

A $\F_1$-scheme $X$ is \emph{of finite type}, if it has a finite covering by affine schemes $U_i=\spec(A_i)$ such that each $A_i$ is a finitely generated monoid.
For a ring $R$, we write $X_R$ for the $R$-base-change of $X$, so $X_R=X_\Z\otimes_\Z R$.

For a monoid $A$ we let $A\otimes\Z$ be the monoidal ring $\Z[A]$.
This defines a functor from monoids to rings which is left adjoint to the forgetful functor that sends a ring $R$ to the multiplicative monoid $(R,\times)$.
This construction is compatible with gluing, so one gets a functor $X\mapsto X_\Z$ from $\F_1$-schemes to $\Z$-schemes.
In \cite{zetaF1} we have shown that 
$X$ is of finite type if and only if $X_\Z$ is a $\Z$-scheme of finite type.

We say that the monoid $A$ is \emph{integral}, if it has the {cancellation property}, i.e., if $ab=ac$ implies $b=c$ in $A$.
This is equivalent to saying that $A$ injects into its quotient group or $A$ is a submonoid of a group.

By a \emph{module} of a monoid $A$ we mean a set $M$ together with a map $A\times M\to M$; $(a,m)\mapsto am$ with $1m=m$ and $(ab)m=a(bm)$.
A \emph{stationary point} of a module is a point $m\in M$ with $am=m$ for every $a\in A$.
A \emph{pointed module} is a pair $(M,m_0)$ consisting of an $A$-module $M$ and a stationary point $m_0\in M$.

\section{Flatness}
Recall the tensor product of two modules $M,N$ of $A$:
$$
M\otimes N\= M\otimes_A N\= M\times N/\sim,
$$
where $\sim$ is the equivalence relation generated by $(am,n)\sim (m,an)$ for every $a\in A, m\in M, n\in N$.
The class of $(m,n)$ is written as $m\otimes n$. 
The tensor product $M\otimes N$ becomes a module via $a(m\otimes n)=(am)\otimes n$.
For example, the module $A\otimes M$ is isomorphic to $M$.

Let now $(M,m_0)$ and $(N,n_0)$ be two pointed modules of $A$, then $(M\otimes N,m_0\otimes n_0)$ is a pointed module, called the pointed tensor product.

The category $\Mod_0(A)$ of pointed modules and pointed morphisms has a terminal and initial object $0$, so it makes sense to speak of kernels and cokernels.
It is easy to see that every morphism $f$ in $\Mod_0(A)$ possesses both.
One defines the \emph{image} of $f$ as $\im(f)=\ker(\coker(f))$ and the \emph{coimage} as $\coim(f)=\coker(\ker(f))$.

A morphism is called \emph{strong}, if the natural map from $\coim(f)$ to $\im(f)$ is an isomorphism.
Kernels and cokernels are strong.
If $A\stackrel{f}{\to}B\stackrel{g}{\to}C$ is given with $g$ being strong and $gf=0$, then the induced map $\coker(f)\to C$ is strong. Likewise, if $f$ is strong and $gf=0$, then the induced map $A\to\ker g$ is strong.
A map is strong if and only if it can be written as a cokernel followed by a kernel.

The usual notion of exact sequences applies, and we say that a sequence of morphisms is \emph{strong exact} if it is exact and all morphisms in the sequence are strong.

A module $F\in\Mod_0(A)$ is called \emph{flat}, if the functor $X\mapsto F\otimes X$ is strong-exact, i.e., if for every strong exact sequence
$$
0\to M\to N\to P\to 0
$$
the induced sequence
$$
0\to F\otimes M\to F\otimes N\to F\otimes P \to 0
$$
is strong exact as well.

It is easy to see that a pointed module $F$ is flat if and only if for every injection $M\hookrightarrow N$ of pointed modules the map $F\otimes M\to F\otimes N$ is an injection.

\emph{Examples.} 
If $A$ is a group, then every module is flat.
Let $S$ be a submonoid of $A$.
Then the localization $S^{-1}A$ is a flat $A$-module. The direct sum $G\oplus F$ of two flat modules is flat. 
Finally, consider the free monoid in one generator $C_+=\{1,\tau,\tau^2,\dots\}$, then an $A$-module $M$ is flat if and only if $\tau m=\tau m'$ implies $m=m'$ for all $m,m'\in M$. This is equivalent to saying that $M$ is a 
$C_+$-submodule of a module of the quotient group $C_\infty=\tau^\Z$ of $C_+$.
The same characterization holds for every integral monoid.

A morphism $\ph :A\to B$ of monoids is called flat if $B$ is flat as an $A$-module.
A morphism of $\F_1$-schemes $f :X\to Y$ is called flat if for every $x\in X$ the morphism of monoids $f^\# : \CO_{Y,f(x)}\to \CO_{X,x}$ is flat.

The following is straightforward.
\begin{itemize}
\item A morphism of monoids $\ph :A\to B$ is flat if and only if the induced morphism of $\F_1$-schemes $\spec B\to\spec A$ is flat.
\item The composition of flat morphisms is flat.
\item The base change of a flat morphism by an arbitrary morphism is flat.
\end{itemize}

{\bf Remark.}
It is easy to see that if $\Z[F]$ is flat as $\Z[A]$-module, then $F$ is flat as $A$-module.
The converse is already false if $A$ is a group.
As an example let $k$ be a field and let $A$ be the group of all matrices of the form $\left(\begin{array}{cc}1 & x \\0 & 1\end{array}\right)$ where $x\in k$.
Let $A$ act on $k^2$ in the usual way and trivially on $k$.
Consider the exact sequence of $\Z[A]$-modules,
$$
\begin{diagram}
\node{0}\arrow{e}
	\node{k}\arrow{e,t}{\al}
		\node{k^2}\arrow{e,t}{\beta}
			\node{k}\arrow{e}
				\node{0,}
\end{diagram}
$$
where $\al(x)=\left(\begin{array}{c}x \\0\end{array}\right)$, $\beta\left(\begin{array}{c}x \\y\end{array}\right)=y$.
Let $F=\{ 1\}$ the trivial $A$-module, then for every $\Z[A]$-module $M$ one has $M\otimes_{\Z[A]}\Z[F]=H_0(A,M)$.
Note that $H_0(A,k)=k$ and that $H_0(\al)=0$, so it is not injective, hence $\Z[F]$ is not flat.

\section{Algebraic extensions}
Let $A$ be a submonoid of $B$.
An element $b\in B$ is called \emph{algebraic over} $A$, if there exists $n\in\N$ with $b^n\in A$.
The extension $B/A$ is called \emph{algebraic},  if every $b\in B$ is {algebraic} over $A$.
An algebraic extension $B/A$ is called \emph{strictly algebraic}, if for every $a\in a$ the equation $x^n=a$ has at most $n$ solutions in $B$.

If $B/A$ is algebraic, then $\Z[B]/\Z[A]$ is an algebraic ring extension, but the converse is wrong in general, as the following example shows: Let $A=\F_1$ and $B$ be the set of two elements, $1$ and $b$ with $b^2=b$.

A monoid $A$ is called \emph{algebraically closed}, if every equation of the form $x^n=a$ with $a\in A$ has a solution in $A$.
Every monoid $A$ can be embedded into an algebraically closed one, and if $A$ is a group, then there exists a smallest such embedding, called the \emph{algebraic closure} of $A$.
 For example, the algebraic closure $\bar\F_1$ of $\F_1$ is the group $\mu_\infty$ of all roots of unity, which is isomorphic to $\Q/\Z$.

\section{Etale morphisms}
Recall that a homomorphism $\ph\colon A\to B$ of monoids is called a \emph{local} homomorphism, if $\ph^{-1}(B^\times)=A^\times$ (every $\ph$ satisfies ``$\supset$'').
For a monoid $A$ let $m_A=A\setminus A^\times$ be its maximal ideal.
It is easy to see that
a homomorphism $\ph\colon A\to B$ is local if and only if $\ph(m_A)\subset m_B$.

A local homomorphism $\ph\colon A\to B$ is called \emph{unramified} if
\begin{itemize}
\item $\ph(m_A)B\= m_B$ and
\item $\ph$ injects $A^\times$ into $B^\times$ and 
$B /\ph(A)$ 
 is a finite strictly algebraic extension. 
\end{itemize}
Note that if $\ph$ is unramified, then so are all localizations $\ph_\p : A_{\ph^{-1}(\p)}\to B_\p$ for $\p\in\spec B$.

A morphism $f:X\to Y$ of $\F_1$-schemes is called unramified, if for every $x\in X$ the local morphism $f^\# : \CO_{Y,f(x)}\to \CO_{X,x}$ is unramified.
 
A morphism $f:X\to Y$ of $\F_1$-schemes is called \emph{locally of finite type}, if every point in $Y$ has an open affine neighborhood $V=\spec A$ such that $f^{-1}(V)$ is a union of open affines $\spec B_i$ with $B_i$ finitely generated as a monoid over $A$.
The morphism is \emph{of finite type} if for every point in $Y$ the number of $B_i$ can be chosen finite.
The morphism is called \emph{finite}, if every $y\in Y$ has an open affine neighborhood $V=\spec A$ such that $f^{-1}(V)$ is affine, equal to $\spec B$, where $B$ is finitely generated as $A$-module.

A morphism $f:X\to Y$ of finite type is called \emph{\'etale}, if $f$ is flat and unramified.
It is called an \emph{\'etale covering}, if it is also finite.

\begin{proposition}
The \'etale coverings of $\spec\F_1$ are the morphisms of the form $\spec A\to\spec\F_1$, where $A$ is a finite cyclic group. The scheme $\spec\bar\F_1$ has no non-trivial \'etale coverings.
\end{proposition}

\prf
Clear.\qed

A connected scheme over $\F_1$, which has only the trivial \'etale covering, is called \emph{simply connected}.

\begin{proposition}
The schemes $\spec\bar\F_1$,  $\spec C_+\times_{\F_1}\bar\F_1$ and $\P^1_{\bar\F_1}$ are simply connected.
\end{proposition}

\prf
The first has been dealt with.
For the second, let $A=\mu_\infty\times C_+$.
Then $\spec A=\spec C_+\times_{\F_1}\bar\F_1$.
Let $f:X\to\spec A$ be an \'etale covering.
As $f$ is finite, $X$ is affine, say $X=\spec B$.
Let $\ph:A\to B$ denote the corresponding morphism of monoids.
The space $\spec A$ consists of two points, the generic point $\eta_A$ and the closed point $c_A$.
Likewise, let $\eta_B, c_B$ denote the generic and closed points of $\spec B$.
One has $f(\eta_B)=\eta_A$.
We will show that $f(c_B)=c_A$.
Assume the contrary.
Then $\ph^{-1}(m_B)$ is empty, hence $\ph$ maps $A$ to the unit group $B^\times$.
The localization at the closed point $c_B$ then maps $\mu_\infty\times C_\infty$ to $B^\times$ and is unramified, hence injective.
But as $C_+\to C_\infty$ is not finite, neither can $\ph$ be finite, a contradiction.
So we conclude $f(c_B)=c_A$, and so the corresponding localization, which is $\ph$ itself, is unramified.
Let $s=\ph(\tau)$, where $\tau$ is the generator of $C_+$.
Then $\ph(m_A)B=m_B$ implies $m_B=sB$, and so $B=B^\times\dotcup sB$ (disjoint union).
Also, $B^\times$ is an algebraic extension of $A^\times\cong\mu_\infty$, hence equals $\ph(A^\times)$.
As $B$ is finitely generated and flat as $A$-module, there are $b_1,\dots ,b_r\in B$ with 
$$
sB\= B^\times s^\N \dotcup B^\times s^\N b_1\dotcup\dots\dotcup B^\times s^\N b_r.
$$
If we assume $r>0$, then $b_1$ is algebraic over $\ph(A)=B^\times\dotcup B^\times s^\N$, so let $N$ be the smallest number in $\N$ such that $b_1^N\in\ph(A)$.
Then $b_1^N\notin B^\times\cong \mu_\infty$, because, as the extension is strictly algebraic, then $b_1$ would be in $B^\times$ already.
So $b_1^N\in B^\times s^\N$.
As the group $B^\times$ is divisible, we can replace $b_1$ with a $B^\times$ multiple to get $b_1^N=s^M$ for some $M\in\N$.
Then $b_1\notin B^\times s^\N b_1$, as $b_1=b^*s^kb_1$ leads to $s^M=b_1^N=(b^*)^Ns^{kN+M}$ which contradicts the injectivity of $\ph$.
But then $b_1$ must be in one of the other $B^\times s^\N$-orbits,which contradicts the disjointness of these orbits.
We conclude $r=0$, i.e. $B=B^\times\dotcup B^\times s^\N\cong A$ as claimed.
The assertion for $\P_{\bar\F_1}^1$ is an easy consequence.
\qed

\section{Toric varieties}
Recall a \emph{toric variety} is an irreducible variety $V$ over $\C$ together with an algebraic action of the $r$-dimensional torus $\GL_1^r$, such that $V$ contains an open orbit.

As toric varieties can be constructed via lattices it follows that every toric variety is the lift $X_\C$ of an $\F_1$-scheme $X$.
For integral schemes of finite type there is a converse direction given in the following theorem, which shows that integral $\F_1$-schemes of finite type are essentially the same as toric varieties.

\begin{theorem}
Let $X$ be a connected integral $\F_1$-scheme of finite type.
Then every irreducible component of $X_\C$ is a toric variety.
The components of $X_\C$ are mutually isomorphic as toric varieties.
\end{theorem}

\prf
Let $U=\spec A$ be an open affine subset of $X$.
Let $\eta$ be the generic point of $X$, then the localization $G=A_\eta$ is the quotient group of $A$.
At the same time, $G$ is the stalk $\CO_{X,\eta}$, so $G$ does not depend on the choice of $U$ up to canonical isomorphism.
Let $\ph: A\to G$ be the quotient map, which is injective as $X$ is integral.
The $\C$-algebra homomorphism,
\begin{eqnarray*}
\C[A] &\to & \C[G]\otimes\C[A]\\
a &\mapsto& \ph(a)\otimes a
\end{eqnarray*}
defines an action of the algebraic group $\CG=\spec\C[G]$ on $\spec\C[A]$.
Since this is compatible with the restriction maps of the structure sheaf, we get an algebraic action of the group scheme $\CG$ on $X_\C$.
As $X$ is integral, $\CG=\spec\C[G] =\spec \C[A_\eta]$ also is an open subset $V_\C$ of $X_\C$, and for $U_\C=\spec\C[A]$ the map
$$
\CO(U_\C)\=\C[A]\ 
\begin{array}{c}\ph\\ \longrightarrow\\ {}\end{array}\ \C[G]\=\CO(V_\C)
$$
is the restriction map of the structure sheaf $\CO$ of $X_\C$.
The map $\C[A]\to \C[G]$ is injective and $\C[G]$ has zero Jacobson radical, so it follows that $V_\C$ is dense in $X_\C$, so in particular it meets every irreducible component.
The group $G$ is a finitely generated abelian group, so $G\equiv \Z^r\times F$ for a finite abelian group $F$.
Hence $\CG\equiv \GL_1^r\times F$ as a group-scheme.
As $\CG$ meets every  component of $X_\C$, the latter are permuted by $F$. Whence the claim.
\qed

To formulate the next result, we will briefly recall the standard construction of toric varieties, see \cite{Fulton}.
Let $N$ be a \emph{lattice}, i.e., a group isomorphic to $\Z^n$ for some $n$. 
A \emph{fan} $\Delta$ in $N$ is a finite collection of \emph{proper convex rational polyhedral cones} $\sigma$ in the real vector space $N_\R=N\otimes\R$, such that every face of a cone in $\Delta$ is in $\Delta$ and the intersection of two cones in $\Delta$ is a face of each.
(Here zero is considered a face of every cone.)
We explain the notation further: A \emph{convex cone} is a convex subset $\sigma$ of $N_\R$ with $\R_{\ge 0}\sigma =\sigma$, it is \emph{polyhedral}, if it is finitely generated and \emph{rational}, if the generators lie in the lattice $N$.
Finally, a cone is called \emph{proper} if it does not contain a non-zero sub vector space of $N_\R$.

Let a fan $\Delta$ be given.
Let $M=\Hom(N,\Z)$ be the dual lattice.
for a cone $\sigma\in \Delta$ the \emph{dual cone} 
$\breve\sigma$ is the cone in the dual space $M_\R$ 
consisting of all $\al\in M_\R$ such that $\al(\sigma)\ge 0$.
This defines a monoid $A_\sigma=\breve\sigma\cap M$.
Set $U_\sigma=\spec(\C[A_\sigma])$.
If $\tau$ is a face of $\sigma$, then $A_\tau\supset A_\sigma$, and this inclusion gives rise to an open embedding $U_\tau\hookrightarrow U_\sigma$.
Along these embeddings we glue the affine varieties $U_\sigma$ to obtain a variety $X_\Delta$ over $\C$, which has a given $\F_1$-structure.
Then $X_\Delta$ is a toric variety, the torus being $U_{0}\cong \GL_1^n$.
Every toric variety is given in this way.

\begin{lemma}
Let $B$ be a submonoid of the monoid $A$ of finite index.
Then the map $\psi:\spec A\to\spec B$ defined by $\psi(\p)=\p\cap B$ is a bijection.
\end{lemma}

\prf
Let $N\in\N$ be such that $a^N\in B$ for every $a\in A$.
To see injectivity, let $\psi(\p)=\psi(\q)$ and let $a\in\p$.
Then $a^N\in\q$ and so $a\in\q$ as $\q$ is a prime ideal.
This shows $\p\subset\q$ and by symmetry we get equality.
For surjectivity, let $\p_B\in\spec B$ and let $\p=\{ a\in A: a^N\in\p_B\}$.
Then $\psi(\p)=\p_B$.
\qed

\begin{proposition}
Suppose that $\Delta$ is a fan in a lattice of dimension $n$.
For $j=0,\dots,n$ let $f_j$ be the number of cones in $\Delta$ of dimension $j$.
Set 
$$
c_j\=\sum_{k=j}^n f_{n-k}(-1)^{k+j}\left(\begin{array}{c}k \\j\end{array}\right).
$$
Let $X$ be the corresponding toric variety, then the $\F_1$-zeta function of $X$ equals
$$
\zeta_X(s)\= s^{c_0}(s-1)^{c_1}\cdots (s-n)^{c_n}.
$$
\end{proposition}

\prf
Let $\sigma\in \Delta$ be a cone of dimension $k$.
Let $F$ be a face of $\breve\sigma$.
Let $\p_F=A_\sigma\setminus F$.
Then $\p_F$ is a non-empty prime ideal in $A_\sigma$.
The map $F\mapsto \p_F$ is a bijection between the set of all faces of $\breve\sigma$ and the set of non-empty prime ideals of $A_\sigma$.
The set $S_\p= A\setminus\p$ equals $M\cap F.$
The quotient group $\Quot(S_\p)$ is isomorphic to $\Z^f$, where $f$ is the dimension of $F$.
There is a bijection between the set of faces of $\sigma$ and the set of faces of $\breve\sigma$ mapping a face $\tau$ to the face $F$ of all $\al\in\breve\sigma$ with $\al(\tau)=0$.
The dimension of $F$ then equals $n-\dim(\tau)$.
So let $f_j^\sigma$ denote the number of faces of $\sigma$ of dimension $j$.
Then the zeta polynomial of $X_\sigma$ equals
$$
N_\sigma(x)\=\sum_{k=0}^n f_k^\sigma (x-1)^{n-k}.
$$
Let $N_{\Delta}$ be the zeta polynomial of $X_{\Delta}$.
We get
\begin{eqnarray*}
N_\Delta(x) &=& \sum_{k=0}^nf_k (x-1)^{n-k}\\
&=& \sum_{k=0}^n f_k \sum_{j=0}^{n-k} \left(\begin{array}{c}n-k \\j\end{array}\right) x^j(-1)^{n-k-j}\\
&=& \sum_{k=0}^n f_{n-k}\sum_{j=0}^k \left(\begin{array}{c}k \\j\end{array}\right)x^j(-1)^{k-j}\\
&=& \sum_{j=0}^n x^j\sum_{k=j}^n f_{n-k}\left(\begin{array}{c}k \\j\end{array}\right) (-1)^{k-j}.
\end{eqnarray*}
This implies the claim.
\qed

\section{Valuations}
On the infinite cyclic monoid $C_+=\{ 1,\tau,\tau^2,\dots\}$ we have a natural linear order given by $\tau^k\le\tau^l\Leftrightarrow k\le l$.
Let $\ph,\psi$ be two monoid morphisms from a monoid $A$ to $C_+$.
Then define $\ph\le\psi\Leftrightarrow \ph(a)\le\psi(a)\ \forall a\in A$.
A \emph{valuation} on $A$ is a non-trivial homomorphism $v:A\to C_+$ which is minimal with respect to the order $\le$ among all non-trivial homomorphisms from $A$ to $C_+$.
Let $V(A)$ denote the set of valuations on $A$.

\begin{lemma}
Let
$$
\begin{diagram}
\node{1}\arrow{e}
	\node{A}\arrow{e}
		\node{B}\arrow{e,t}{\ph}
			\node{F}\arrow{e}
				\node{1}
\end{diagram}
$$
be an exact sequence of monoids, where $F$ is a finite abelian group.
Then for every valuation $v\in V(A)$ there exists a unique valuation $w$ on $B$ and $k\in\N$ such that
$$
w|_A\= v^k.
$$
Mapping $v$ to $w$ sets up a bijection from $V(A)$ to $V(B)$.
\end{lemma}

\prf
Let $F'$ be a subgroup of $F$ and let $B'$ be the preimage of $F'$ under $\ph$.
We get two exact sequences
$$
\begin{diagram}
\node{1}\arrow{e}
	\node{A}\arrow{e}
		\node{B'}\arrow{e}
			\node{F'}\arrow{e}
				\node{1,}
\end{diagram}
$$
and
$$
\begin{diagram}
\node{1}\arrow{e}
	\node{B'}\arrow{e}
		\node{B}\arrow{e}
			\node{F/F'}\arrow{e}
				\node{1.}
\end{diagram}
$$
Assume we have proven the lemma for each of these two sequences, then it follows for the original one.
In this way we reduce the proof to the case when $F$ is a finite cyclic group.
We first show existence of $w$ for given $v$.
For this let $f_0$ be a generator of $F$ and let $l$ be its order.
Choose a $b_0$ in the preimage $\ph^{-1}(f_0)$. Then $b_0^l\in A$, and $v(b_0^l)=\tau^n$ for some $n\ge 0$.
If $n=0$, then set $k=1$ and define $w: B\to C_+$ by $w(b_0^ja)=v(a)$ for $a\in A$ and $j\ge 0$.
If $n>0$, then set $k=l/{\rm gcd(l,n)}$ and let $w: B\to C_+$ be defined by $w(b_0^ja)=\tau^j v(a)^k$.
This shows existence of the extension $w$.
\qed

\section{Cohomology}
Cohomology is not defined over $\F_1$.
I am grateful to Ofer Gabber for bringing the following example to my attention. 
Let $X$ be the topological space consisting of three points $\eta, X_+,x_-$.
The open sets besides the trivial ones are $U=\{ \eta\},  U_+=\{\eta,x_+\}, U_-=\{\eta, x_-\}$.
Let $A$ be a subgroup of the abelian group $B$ and let $C=B/A$. 
Let $\CF$ be the sheaf of abelian groups on $X$ with $\CF(U_\pm)=A$ and $\CF(U)=B$ and the restriction being the inclusion.
Let $\CG$ be the constant sheaf $B$ and let $\CH$ be the quotient sheaf $\CG/\CF$.
As $\CG$ is flabby, the long cohomology sequence terminates and looks like this:
$$
0\to H^0(\CF)\to H^0(\CG)\to H^0(\CH)\to H^1(\CF)\to 0 
$$
In concrete terms this is
$$
0\to A\to B\to C \times C\to (C\times C)/\Delta\to 0,
$$
where $\Delta$ means the diagonal in $C\times C$.
Let $f: X\to X$ be the homeomorphism with $f(x_+)=x_-$, $f(x_-)=x_+$, and $f(\eta)=\eta$.
There is a natural isomorphism $f_*\CF\cong\CF$ and for the other sheaves as well.
On the global sections of $\CF$ and $\CG$ this induces the trivial map, whereas on $H^0(\CH)$ it induces the flip $(a,b)\mapsto (b,a)$, which on $H^1(\CF)$ amounts to the same as the inversion $a\mapsto -a$.
The naturality of these isomorphisms means that if the sheaves and the cohomology groups  are defined over $\F_1$, then so must be the flip.
This, however, is not the case, as for a set $S$ the inversion on the abelian group $\Z[S]$ is not induced by a self-map of $S$.

Even more convincing is the fact that in this example there are different injective resolutions which produce different cohomology groups.

{\small Mathematisches Institut\\
Auf der Morgenstelle 10\\
72076 T\"ubingen\\
Germany\\
\tt deitmar@uni-tuebingen.de}

\end{document}